\documentclass[11pt]{amsart}

\usepackage{amsmath,amsthm,amssymb,amsxtra,amsfonts,amsbsy,mathrsfs}
\usepackage{verbatim}
\usepackage[all]{xy}
\newtheorem{theorem}{Theorem}[section]
\newtheorem{corollary}[theorem]{Corollary}
\newtheorem{lemma}[theorem]{Lemma}
\newtheorem{proposition}[theorem]{Proposition}
\newtheorem{definition}[theorem]{Definition}

\newtheorem{remark}[theorem]{Remark}
\newtheorem{example}[theorem]{Example}

\newcommand{\Z}{\mathbb{Z}}

\newcommand{\N}{\mathbb{N}}

\newcommand{\til}[1]{\tilde{#1}}
\newcommand{\mf}[1]{\mathfrak{#1}}
\newcommand{\OO}{\mathcal{O}}
\newcommand{\kVG}{k[V]^G}
\newcommand{\Vv}{\mathbb{V}^\vee}

\DeclareMathOperator{\Spec}{Spec}

\DeclareMathOperator{\Hom}{Hom}

\DeclareMathOperator{\Tor}{Tor}

\DeclareMathOperator{\Aut}{Aut}

\DeclareMathOperator{\Sym}{Sym}
\DeclareMathOperator{\Gal}{Gal}

\SelectTips{xy}{12}

\makeatletter
\@namedef{subjclassname@2010}{ \textup{2010} Mathematics Subject 
Classification} \makeatother

\begin{document}

\keywords{Chevalley-Shephard-Todd, pseudo-reflection, linearly reductive, tame stacks}
\subjclass[2010]{14A20, 14L15.}

\title[CST for Finite Linearly Reductive Group Schemes]{The Chevalley-Shephard-Todd Theorem for Finite Linearly 
Reductive Group Schemes}
\author{Matthew Satriano}
\date{}
\maketitle

\vspace{-1cm}

\begin{abstract}
We generalize the classical Chevalley-Shephard-Todd theorem to the case of finite linearly reductive group schemes. 
As an application, we prove that every scheme $X$ which is \'etale locally the quotient of a smooth scheme by a finite 
linearly reductive group scheme is the coarse space of a smooth tame Artin stack (as defined by Abramovich, Olsson, 
and Vistoli) whose stacky structure is supported on the singular locus of $X$. 
\end{abstract}

\vspace{0.5cm}

\section{Introduction}
\label{sec:intro}
Given a field $k$ and an action of a finite (abstract) group $G$ on a $k$-vector space $V$, we obtain a linear action of $G$ 
on the polynomial ring $k[V]$.  A central theme in Invariant Theory is determining when certain nice properties of a ring 
with $G$-action are inherited by its invariants.  In particular, it is natural to ask when $\kVG$ is polynomial.  
If $G$ acts faithfully on $V$, we say $g\in G$ is a 
\emph{pseudo-reflection} (with respect to the action of $G$ on $V$) if $V^g$ is a hyperplane.  The classical 
Chevalley-Shephard-Todd Theorem states
\begin{theorem}[{\cite[\S 5 Thm 4]{bourbaki}}]
\label{thm:classicalcst}
If $G\rightarrow \Aut_k(V)$ is a faithful representation of a finite group and the order of $G$ is not divisible by 
the characteristic of $k$, then $\kVG$ is polynomial if and only if $G$ 
is generated by pseudo-reflections.
\end{theorem}
In this paper we generalize this theorem to the case of finite linearly reductive group schemes.  To do so, we first need 
a notion of pseudo-reflection in this setting.
\begin{definition}
\label{def:prlin}
\emph{Let $k$ be a field and $V$ a finite-dimensional $k$-vector space with a faithful action of a finite linearly reductive group scheme 
$G$ over $\Spec k$.  We say that a subgroup scheme $N$ of $G$ is a \emph{pseudo-reflection} if $V^N$ has codimension 
$1$ in $V$.  We define the \emph{subgroup\ scheme\ generated\ by\ pseudo-reflections} to be the intersection of the subgroup schemes 
which contain all of the pseudo-reflections of $G$.  We say $G$ is \emph{generated\ by\ pseudo-reflections} if $G$ is the 
subgroup scheme generated by pseudo-reflections.}
\end{definition}
Over algebraically closed fields, Theorem \ref{thm:classicalcst} generalizes to 
\begin{theorem}
\label{thm:algclosed}
Let $k$ be an algebraically closed field and $V$ a finite-dimensional $k$-vector space with a faithful action of a finite 
linearly reductive group scheme $G$ over $\Spec k$.  Then $G$ is generated by pseudo-reflections if and only if $k[V]^G$ 
is polynomial.
\end{theorem}
A more technical version of this theorem holds over fields which are not algebraically closed; however, the ``only if'' direction 
does not hold for finite linearly reductive group schemes in general (see Example \ref{ex:prdescend}).  We instead 
prove the ``only if'' direction for the smaller 
class of stable group schemes, which we now define (see Proposition \ref{prop:idex} for examples).  Over an algebraically 
closed field, the class of stable group schemes coincides with that of finite linearly reductive group schemes.  Recall from \cite[Def 2.9]{tame} that $G$ is called \emph{well-split} if it is isomorphic to a semi-direct product 
$\Delta\rtimes Q$, where $\Delta$ is a finite diagonalizable group scheme and $Q$ is a finite constant tame group scheme; here, 
tame means that the degree is prime to the characteristic.
\begin{definition}
\label{def:id}
\emph{
A group scheme $G$ over a field $k$ is called \emph{stable} if the following two conditions hold:
\begin{enumerate}
\renewcommand{\labelenumi}{(\alph{enumi})}
\item for all finite field extensions $K/k$, every subgroup scheme of $G_K$ descends to a subgroup scheme of $G$
\item there exists a finite Galois extension $K/k$ such that $G_K$ is well-split.
\end{enumerate}
}
\end{definition}
\begin{remark}
\label{rmk:automatic}
\emph{
If $G$ is a finite linearly reductive group scheme over a perfect field $k$, then \cite[Lemma 2.11]{tame} shows that 
condition (b) above is automatically satisfied.
}
\end{remark}
Theorem \ref{thm:algclosed} is then a special case of the following generalization of the Chevalley-Shephard-Todd theorem.  
This is the first main result of this paper.
\begin{theorem}
\label{thm:newcst}
Let $k$ be a field and $V$ a finite-dimensional $k$-vector space with a faithful action of a finite linearly 
reductive group scheme $G$ over $\Spec k$.  If $G$ is generated by pseudo-reflections, then $k[V]^G$ is polynomial.  
The converse holds if 
$G$ is stable.
\end{theorem}
We also prove a version of this theorem for an action of a finite linearly reductive group scheme on a smooth scheme.
\begin{definition}
\label{def:pr}
\emph{Given a smooth affine scheme $U$ over $\Spec k$ with a faithful action of a finite linearly reductive group scheme $G$ which fixes a 
field-valued point $x\in U(K)$, we say a subgroup scheme $N$ of $G$ is a \emph{pseudo-reflection} \emph{at} $x$ if 
$N_K$ is a pseudo-reflection with respect to the induced action of $G_K$ on the cotangent space 
at $x$.  We define what it means for $G$ to be generated by pseudo-reflections at $x$ in the same manner as in Definition \ref{def:prlin}.
}
\end{definition}
Theorem \ref{thm:newcst} then has the following corollary.
\begin{corollary}
\label{cor:newcst}
Let $k$ be a field, let $U$ be a smooth affine $k$-scheme with a faithful action by a finite linearly 
reductive group scheme $G$ over $\Spec k$.  Let $x\in U(K)$, where $K/k$ is a finite separable field extension, and 
suppose $x$ is fixed by $G$.  If $G$ is 
generated by pseudo-reflections at $x$, then $U/G$ is smooth at the image of $x$.  The converse holds if 
$G$ is stable.
\end{corollary}
The second main result of this paper is
\begin{theorem}
\label{thm:torsor}
Let $k$ be a 
field and let $U$ be a smooth affine $k$-scheme with a faithful action by a stable 
group scheme $G$ over $\Spec k$.  Suppose $K/k$ is a finite separable field extension and $G$ fixes a point $x\in U(K)$.  
Let $M=U/G$, let $M^0$ be the smooth 
locus of $M$, and let $U^0=U\times_M M^0$.  If $G$ has no pseudo-reflections at $x$, then after possibly shrinking $M$ 
to a smaller Zariski neighborhood of 
the image of $x$, we have that $U^0$ is a $G$-torsor over $M^0$.
\end{theorem}
We remark that in the classical case, Theorem \ref{thm:torsor} follows directly from Corollary \ref{cor:newcst} and 
the purity of the branch locus theorem \cite[X.3.1]{sga1}.  For us, however, a little more work is needed since $G$ is 
not necessarily \'etale.
\\
\\
As an application of Theorem \ref{thm:torsor}, we 
generalize the 
well-known result (see for example \cite[2.9]{int} or \cite[Rmk 4.9]{fant}) 
that schemes with quotient singularities prime to the characteristic are coarse spaces of smooth 
Deligne-Mumford stacks.  We say a scheme has 
\emph{linearly\ reductive\ singularities} if it is \'etale locally 
the quotient of a smooth scheme by a finite linearly reductive group scheme.  We 
show that every such scheme $M$ is the coarse space of a smooth tame Artin stack 
(in the sense of \cite{tame}) whose stacky structure is supported at the 
singular locus of $M$.  More precisely,
\begin{theorem}
\label{prop:cslrs}
Let $k$ be a perfect field and $M$ a $k$-scheme with linearly reductive singularities.  Then it is the coarse space of 
a smooth tame stack $\mf{X}$ over $k$ such that $f^0$ in the diagram 
\[
\xymatrix{
\mf{X}^0\ar[r]^{j^0}\ar[d]_{f^0} & \mf{X}\ar[d]^f \\
M^0\ar[r]_j & M
}
\]
is an isomorphism, where $j$ is the inclusion of the smooth locus of $M$ 
and $\mf{X}^0=M^0\times_{M}\mf{X}$.
\end{theorem}
%
This paper is organized as follows.  In Section \ref{sec:polynomial}, we prove the ``if'' direction of 
Theorem \ref{thm:newcst} and reduce the proof of the ``only if'' direction to the special case of 
Theorem \ref{thm:torsor} in which $U=\Vv(V)$ for some $k$-vector space $V$ with $G$-action (see the Notation section 
below).  This special case is proved in Section \ref{sec:polytorsor}.  The key input for the proof is a result of 
Iwanari \cite[Thm 3.3]{mfr} which we reinterpret in the 
language of pseudo-reflections.  We finish the section by proving Corollary \ref{cor:newcst}.  
In Section \ref{sec:torsor}, we use Corollary \ref{cor:newcst} to complete the proof of Theorem \ref{thm:torsor}.  
In Section \ref{sec:lrs}, we prove Theorem \ref{prop:cslrs}.\\
\\
$\textbf{Acknowledgments.}$  I would like to thank Dustin Cartwright, Ishai Dan-Cohen, Anton Geraschenko, and David Rydh 
for many helpful conversations.  I am of course indebted to my advisor, Martin Olsson, for his suggestions and 
help in editing this paper.  I thank Dan Edidin for suggesting that I write up these results in this 
stand-alone paper rather than include them in \cite{dR}. Lastly, I would like to thank the referee for his helpful comments.\\
\\
$\textbf{Notation.}$  Throughout this paper, $k$ is a field and $S=\Spec k$.  If $V$ is a $k$-vector space with 
an action of a group scheme $G$, then 
we denote by $\Vv(V)$, or simply $\Vv$ if $V$ is understood, 
the scheme $\Spec k[V]$ whose $G$-action is given by the dual representation on functor points.  Said another way, if
$G=\Spec A$ is affine and its action on $V$ is given by the co-action map $\sigma:V\longrightarrow V\otimes_k A$, then 
the co-action map $k[V]\longrightarrow k[V]\otimes_k A$ defining the $G$-action on $\Vv$ is given by 
$\sum a_iv_i\mapsto \sum a_i\sigma(v_i)$.\\
\\
All Artin stacks $\mf{X}$ in this paper are assumed to have finite diagonal so that if $\mf{X}$ is locally of finite presentation, it has a 
coarse space by \cite[Thm 1.1]{kmc} (\emph{c.f.} \cite{km}).  Given a locally finitely presented scheme $U$ with an action of a finite flat group scheme $G$, we 
denote by $U/G$ the coarse space of the stack $[U/G]$.\\
\\
If $R$ is a ring and $\mathcal{I}$ an ideal of $R$, then we denote by $V(\mathcal{I})$ the closed subscheme of 
$\Spec R$ defined by $\mathcal{I}$.

\section{Linear Actions on Polynomial Rings}
\label{sec:polynomial}
\subsection{The ``if'' direction of Theorem \ref{thm:newcst}}
Our goal in this subsection is to prove the ``if'' direction of Theorem \ref{thm:newcst}.  
We begin with examples of stable group schemes and 
with some basic results about the subgroup scheme generated by pseudo-reflections.
\begin{lemma}
\label{l:wsoverZ}
Suppose $k$ is perfect and $G$ is a finite linearly reductive group scheme over $S$.  If the identity component $\Delta$ of $G$ is diagonalizable 
and $G/\Delta$ is constant, then there exists a finite linearly reductive group scheme $\til{G}$ over $\Z$ such that $\til{G}_k=G$.  
If $H$ is a closed subgroup scheme of $G$, then there exists a closed subgroup scheme $\til{H}$ of $\til{G}$ whose pullback to $k$ is $H$.  
If $H$ is normal in $G$, then $\til{H}$ is normal in $\til{G}$.
\end{lemma}
\begin{proof}
Let $Q=G/\Delta$.  Since $k$ is perfect, the connected-\'etale sequence
\[
1\longrightarrow \Delta\longrightarrow G\longrightarrow Q\longrightarrow 1
\]
is functorially split (see \cite[3.7 (IV)]{tate}).  Since $\Delta$ is diagonalizable, 
it is of the form $\Spec k[A]$, where $A$ is a finitely generated abelian group.  Note that as a scheme $G=\Delta\times_k Q$ and that 
its group scheme structure is given by a homomorphism
\[
\epsilon:Q\longrightarrow\mathcal{A}ut(\Delta)=\Aut(A).
\]
We can therefore let $\til{G}=\Spec\Z[A]\times_\Z Q$ with group scheme structure induced by $\epsilon$.\\
\\
Now let $H$ be a closed subgroup scheme of $G$.  Letting $\Delta'=H \cap\Delta$ and $Q'=H/\Delta'$, we have a commutative diagram
\[
\xymatrix{
1\ar[r] & \Delta\ar[r] & G\ar[r] & Q\ar[r] & 1\\
1\ar[r] & \Delta' \ar[r]\ar[u]^\varphi & H\ar[r]\ar[u] & Q'\ar[r]\ar[u]_\psi & 1
}
\]
with exact rows.  Since $\Delta$ is connected, we see $\Delta'$ is the connected component of the identity of $H$.  
Therefore, the bottom row of the above diagram is the connected-\'etale sequence of $H$, and so
\[
H=\Delta'\rtimes Q'
\]
as $k$ is perfect.  Since $\Delta'$ is diagonalizable and $Q'$ is constant, we can define $\til{H}$ in the same way we defined $\til{G}$.\\
\\
We now show that $\til{H}$ is a closed subgroup scheme of $\til{G}$.  Let $*$ denote the action of $Q$ (resp. $Q'$) on 
$\Delta$ (resp. $\Delta'$).  Since the splitting of the connected-\'etale sequence of a finite group scheme over a 
perfect field is functorial, 
we see that for all $q'\in Q'$ and local sections $\delta'$ of $\Delta'$,
\[
\psi(q')*\varphi(\delta')=\varphi(q'*\delta').
\]
We therefore obtain a closed immersion from $\til{H}$ to $\til{G}$ whose pullback to $k$ is the morphism from $H$ to $G$.\\
\\
Lastly, we show that if $H$ is normal in $G$, then $\til{H}$ is normal in $\til{G}$.  Let $\Delta'=\Spec k[A']$, where $A'$ is a finitely-generated abelian 
group.  Showing that $\til{H}$ is normal in $\til{G}$ is equivalent to showing that $Q'$ is normal in $Q$, and for all local sections $\delta\in\Delta$, 
$\delta'\in\Delta'$, $q\in Q$, and $q'\in Q'$, we have
\[
q*(\delta^{-1}\delta'\cdot({q'}^{-1}*\delta))\in\Delta'.
\]
We know that $Q'$ is normal in $Q$ as $H$ is normal in $G$.  To check the latter statement about local sections, note that it can be reformulated 
as follows: for every $q\in Q$ and $q'\in Q$, the homomorphism
\[
A\longrightarrow A\times A'
\]
\[
a\longmapsto(q*(a^{-1}\cdot {q'}^{-1}*a),q*\bar{a})
\]
factors through $A'$; here $\bar{a}$ denotes the image of $a$ under the projection from $A$ to $A'$.  Since this statement makes no reference 
to the base scheme, it can be checked over $k$, where the normality of $H$ in $G$ yields the desired factorization.
\end{proof}
\begin{proposition}
\label{prop:idex}
Let $G$ be a finite group scheme over $S$.  Consider the following conditions:
\begin{enumerate}
\item $G$ is diagonalizable.
\item $G$ is a constant group scheme.
\item $k$ is perfect, the identity component $\Delta$ of $G$ is diagonalizable, and $G/\Delta$ is constant.
\end{enumerate}
If any of the above conditions hold, then $G$ is stable.
\end{proposition}
\begin{proof}
It is clear that finite diagonalizable group schemes and finite constant group schemes are stable, so we consider 
the last case.  Let $Q=G/\Delta$.  Since $k$ is perfect, the connected-\'etale sequence
\[
1\longrightarrow \Delta\longrightarrow G\longrightarrow Q\longrightarrow 1
\]
is functorially split.  
Let $K/k$ be a finite extension and let $H$ be a subgroup scheme of 
$G_K$.  Letting $\Delta'=H \cap\Delta_K$ and $Q'=H/\Delta'$, we have a commutative diagram
\[
\xymatrix{
1\ar[r] & \Delta_K\ar[r] & G_K\ar[r] & Q_K\ar[r] & 1\\
1\ar[r] & \Delta' \ar[r]\ar[u] & H\ar[r]\ar[u] & Q'\ar[r]\ar[u] & 1
}
\]
with exact rows.  Since $\Delta$ is connected and has a 
$k$-point, \cite[4.5.14]{ega4} shows that $\Delta$ is geometrically connected.  In particular, $\Delta_K$ is the 
connected component of the identity of $G_K$, and so $\Delta'$ is the connected component of the identity of $H$.  
Therefore, the bottom row of the above diagram is the connected-\'etale sequence of $H$.  
The proposition then follows from Lemma \ref{l:wsoverZ}.
\end{proof}
\begin{lemma}
\label{l:prdescend}
Let $V$ be a finite-dimensional $k$-vector space with a faithful action of a stable group scheme 
$G$ over $S$, and let $H$ be the subgroup scheme generated by 
pseudo-reflections.  If $K/k$ is an algebraic extension of fields, then a subgroup scheme of $G_K$ is a pseudo-reflection 
if and only if it descends to a pseudo-reflection.  Furthermore, $H_K$ is the 
subgroup scheme of $G_K$ generated by pseudo-reflections.\\
\end{lemma}
\begin{proof}
Note first that if $P$ is a subgroup scheme of $G_K$, then there exists a subgroup scheme $P_0$ of $G$ such that $(P_0)_K=P$.  
If $K/k$ is a finite extension, this follows from the fact that $G$ is stable.  If $K/k$ is an infinite extension, by a standard limit argument, 
there exists a finite extension $L/k$ and a subgroup scheme $P_1$ of $G_L$ such that $(P_1)_K=P$.  We then obtain our desired $P_0$ 
as $L/k$ is a finite extension.  The first claim of the proposition then follows from the fact that 
\[
(V_K)^{N_K}=(V^N)_K
\]
for any subgroup scheme $N$ of $G$.  The second claim follows from the fact that 
if $P'$ and $P''$ are subgroup schemes of $G$, then $P'_K$ contains $P''_K$ if and only if $P'$ contains $P''$.
\end{proof}
We remark that even in characteristic zero, Lemma \ref{l:prdescend} is false for general finite linearly reductive 
group schemes $G$, as the following example shows.  Note that this example also shows that the ``only if'' direction of 
Theorem \ref{thm:newcst} and of Corollary \ref{cor:newcst} is false for general finite linearly reductive group schemes.
\begin{example}
\label{ex:prdescend}
\emph{
Let $k$ be a field contained in $\mathbb{R}$ or let $k=\mathbb{F}_p$ for $p$ congruent to $3$ mod $4$.  Let $K=k(i)$, 
where $i^2=-1$, and let $G$ be the locally constant group scheme over $\Spec k$ whose pullback to $\Spec K$ is 
$\Z/2\times\Z/2$ with the Galois action that switches the two $\Z/2$ factors.  Let $g_1$ and $g_2$ be the generators 
of the two $\Z/2$ factors and consider the action
\[
\rho:G_K\longrightarrow \Aut_K(K^2)
\]
on the $K$-vector space $K^2$ given by 
\[
\rho(g_1):(a,b)\mapsto (-bi,ai)
\]
\[
\rho(g_2):(a,b)\mapsto (bi,-ai).
\]
Then $\rho$ is Galois-equivariant and hence comes from an action of $G$ on $k^2$.  Note that $\Z/2\times1$ and 
$1\times\Z/2$ are both pseudo-reflections of $G_K$, as the subspaces which they fix are $K\cdot(1,i)$ and 
$K\cdot(1,-i)$, respectively.  
Since $G_K$ is not a pseudo-reflection, it follows that there are no Galois-invariant pseudo-reflections of $G_K$, 
and hence, the subgroup scheme generated by pseudo-reflections of $G$ is trivial; the subgroup scheme generated by 
pseudo-reflections of $G_K$, however, is $G_K$.
}
\end{example}
\begin{corollary}
\label{l:normal}
If $V$ is a finite-dimensional $k$-vector space with a faithful action of a stable group scheme 
$G$ over $S$, then the subgroup scheme generated by pseudo-reflections is normal in $G$.
\end{corollary}
\begin{proof}
We denote by $H$ the subgroup scheme generated by pseudo-reflections.  
Let $T$ be an $S$-scheme and let $g\in G(T)$.  We must show the subgroup schemes $H_T$ and $gH_Tg^{-1}$ of $G_T$ are equal.  
To do so, it suffices to check this on stalks and so we can assume $T=\Spec R$, where $R$ is strictly Henselian.  By 
\cite[Lemma 2.17]{tame}, we need only show that these two group schemes are equal over the closed fiber of $T$, so 
we can further assume that $R=K$ is a field.  Since $G$ is finite over $S$, the residue fields of $G$ 
are finite extensions of $k$.  We can therefore assume that $K/k$ is a finite field extension.\\
\\
By Lemma  \ref{l:prdescend}, we know that $H_K$ is the subgroup scheme of $G_K$ generated by pseudo-reflections.  
Note that if $N'$ is a pseudo-reflection 
of $G_K$, then $gN'g^{-1}$ is as well since
\[
V_K^{gN'g^{-1}}=g(V_K^{N'}).
\]
As a result, $gH_Kg^{-1}=H_K$, which completes the proof.
\end{proof}
\begin{lemma}
\label{l:genbypr}
Given a finite-dimensional $k$-vector space $V$ with a faithful action of a finite linearly reductive group scheme 
$G$ over $S$, let $\{N_i\}$ denote the set of pseudo-reflections of $G$ and let $H$ be the subgroup scheme generated by 
pseudo-reflections.  Then
\[
k[V]^H=\bigcap_i k[V]^{N_i}.
\]
\end{lemma}
\begin{proof}
Let $R=\bigcap_i k[V]^{N_i}$.  Consider the functor
\[
F:(k\textrm{-}alg)\longrightarrow (Groups)
\]
\[
A\longmapsto\{g\in G(A)\,\mid\, g(m)=m\textrm{\ for\ all\ } m\in R\otimes_k A\}.
\]
Since each $k[V]^{N_i}$ is finitely generated, we see $R$ is as well.  Let $r_1,\dots,r_n$ be a finite set of generators for $R$.  
We see then that $F$ is representable by the intersection of the stabilizers $G_{r_j}$, and so is 
a closed subgroup scheme of $G$.  Since $F$ contains every pseudo-reflection, we see $H\subset F$.  We therefore 
have the containments
\[
R\subset k[V]^{F}\subset k[V]^H\subset \bigcap_i k[V]^{N_i}
\]
from which the lemma follows.
\end{proof}
If $N$ is any subgroup scheme of $G$, it 
is linearly reductive by \cite[Prop 2.7]{tame}.  It follows that 
\[
V\simeq V^N\oplus V/V^N
\]
as $N$-representations.  If $N$ is a pseudo-reflection, then $\dim_k V/V^N=1$.  Let $v$ be a generator of 
the $1$-dimensional subspace $V/V^N$ and let $\sigma:V\rightarrow V\otimes_k B$ be the coaction map, where 
$N=\Spec B$.  Then via the above isomorphism, $\sigma$ is given by 
\[
V^N\oplus V/V^N\longrightarrow (V^N\otimes_k B)\oplus (V/V^N\otimes_k B)
\]
\[
(w,w')\longmapsto (w\otimes 1, w' \otimes b)
\]
for some $b\in B$.  It follows that there is a $k$-linear map $h:V\rightarrow B$ such that for all $w\in V$, 
\[
\sigma(w)-(w\otimes1)=v\otimes h(w).
\]
If we continue to denote by $\sigma$ the induced coaction map $k[V]\longrightarrow k[V]\otimes_k B$, we see that $h$ extends 
to a $k[V]^N$-module homomorphism $k[V]\longrightarrow k[V]\otimes_k B$, which we continue to denote by $h$, such that for 
all $f\in k[V]$, 
\[
\sigma(f)-(f\otimes1)=(v\otimes1)\cdot h(f).
\]
We are now ready to prove the ``if'' direction of Theorem \ref{thm:newcst}.  Our proof is only a slight variant of the 
proof of the classical Chevalley-Shephard-Todd Theorem presented in \cite{cst}.
\begin{proof}[Proof of ``if'' direction of Theorem \ref{thm:newcst}]
Let $R=k[V]^G$.  
By Lemma \ref{l:genbypr}, we know that the intersection of the $k[V]^N$ is $R$, where $N$ runs 
through the pseudo-reflections of $G$.  
By the proposition on page 225 of \cite{cst}, to show $R$ is polynomial, we need only show that $k[V]$ is a free 
$R$-module.  By graded Nakayama, the projective dimension of $k[V]$ is the smallest integer $i$ such that 
$\Tor^{R}_{i+1}(k,k[V])=0$, where $k$ is viewed as an $R$-module via the augmentation map
\[
\epsilon:\kVG\longrightarrow k[V]\longrightarrow k
\]
sending all positively graded elements to $0$.  We must therefore show $\Tor^{R}_1(k,k[V])=0$.\\
\\
Tensoring the short exact sequence defined by $\epsilon$ with $k[V]$, we obtain a long exact sequence 
\[
0\longrightarrow\Tor^{R}_1(k,k[V])\longrightarrow\ker\epsilon\otimes_{R}k[V]\stackrel{\phi}{\longrightarrow} 
R\otimes_{R} k[V]\stackrel{\epsilon\otimes1}{\longrightarrow} k\otimes_{R}k[V]\longrightarrow 0.
\]
To show $\Tor^{R}_1(k,k[V])=0$, we must prove that $\phi$ is injective.  We in fact show 
\[
\phi\otimes 1:\ker\epsilon\otimes_{R}k[V]\otimes_k C\longrightarrow k[V]\otimes_k C
\]
is injective for all finite-dimensional $k$-algebras $C$.  If this is not the case, then the set  
\[
\{\xi\,\mid\, C\textrm{\ is\ a\ finite-dimensional\ } k\textrm{-algebra},\, 0\neq\xi\in 
\ker\epsilon\otimes_{R}k[V]\otimes_k C,\,(\phi\otimes1)(\xi)=0 \}
\]
is non-empty and we can choose an element $\xi$ of minimal degree, where $\ker\epsilon$ is given its natural grading 
as a submodule of $k[V]$ and the elements of $C$ are defined to be of degree 0.  We begin by showing 
$\xi\in \ker\epsilon\otimes_{R}R\otimes_k C$.  That is, we show $\xi$ is fixed by all pseudo-reflections.\\
\\
Let $N=\Spec B$ be a pseudo-reflection.  Let $\sigma:k[V]\longrightarrow k[V]\otimes B$ be the coaction map.  As explained 
above, we get a $k[V]^N$-module homomorphism $h:k[V]\longrightarrow k[V]\otimes B$.  Note that this morphism has 
degree -1.  Since 
\[
(1\otimes \sigma\otimes 1)(\xi)-\xi\otimes1=(1\otimes h\otimes 1)(\xi)\cdot(1\otimes v\otimes1\otimes1),
\]
the commutativity of 
\[
\xymatrix{
\ker\epsilon\otimes k[V]\otimes B\otimes C \ar[rr]^-{\phi\otimes 1\otimes 1} & & k[V]\otimes B\otimes C \\
\ker\epsilon\otimes k[V]\otimes C\ar[u]_{1\otimes \sigma\otimes 1}\ar[rr]^-{\phi\otimes 1} & & 
k[V]\otimes C\ar[u]^{\sigma\otimes 1}
}
\]
implies 
\[
(\phi\otimes 1\otimes 1)(1\otimes h\otimes 1)(\xi)\cdot(v\otimes1\otimes1)=0.
\]
It follows that $(1\otimes h\otimes 1)(\xi)$ is killed by $\phi\otimes 1\otimes 1$.  Since $h$ has degree -1, our 
assumption on $\xi$ shows that $(1\otimes h\otimes 1)(\xi)=0$.  We therefore have 
$(1\otimes \sigma\otimes 1)(\xi)=\xi\otimes1$, which proves that $\xi$ is $N$-invariant.\\
\\
Since $G$ is linearly reductive, we have a section of the inclusion $k[V]^G\hookrightarrow k[V]$.  We therefore, 
also obtain a section $s$ of the inclusion $j:R\hookrightarrow k[V]$.  Let 
$\psi:\ker\epsilon\otimes_{R} R\longrightarrow R$ be the canonical map, and consider the diagram 
\[
\xymatrix{
\ker\epsilon\otimes k[V]\otimes C \ar[rr]^-{\phi\otimes1}\ar@<+.5ex>[d]_{1\otimes j\otimes1\ } & & 
k[V]\otimes C \ar@<+.5ex>[d]_{j\otimes1\ }\\
\ker\epsilon\otimes R\otimes C \ar[rr]^-{\psi\otimes1}\ar@<+.5ex>[u]_{\ 1\otimes s\otimes 1} & & 
R\otimes C \ar@<+.5ex>[u]_{\ s\otimes1}
}
\]
We see that 
\[
(j\otimes1)(\psi\otimes1)(1\otimes s\otimes1)(\xi)=(\phi\otimes1)(1\otimes j\otimes1)(1\otimes s\otimes1)(\xi)=
(\phi\otimes1)(\xi)=0.
\]
But $j\otimes1$ and $\psi\otimes1$ are injective, so $(1\otimes s\otimes 1)(\xi)=0$.  Since 
$\xi\in \ker\epsilon\otimes_{R}R\otimes_k C$, it follows that $\xi=0$, which is a contradiction.
\end{proof}
\subsection{Reducing the ``only if'' direction of Theorem \ref{thm:newcst} to a case of Theorem \ref{thm:torsor}}
Now that we have proved the ``if'' direction of Theorem \ref{thm:newcst}, we work toward reducing the ``only if'' 
direction to 
the special case of Theorem \ref{thm:torsor} where $U=\Vv$.  
The main step in this reduction is showing that if 
$G$ acts faithfully on $V$, and $H$ denotes the subgroup scheme generated by pseudo-reflections, then the action of 
$G/H$ on $\Vv/H$ has no pseudo-reflections at the origin.  In the classical case, the proof of this statement relies on 
the fact that $G$ has no pseudo-reflections if and only if $\Vv\rightarrow\Vv/G$ is \'etale in codimension one.  As the following 
example illustrates, this relation between pseudo-reflections and ramification no longer holds in our case.
\begin{example}
\emph{
Let $k$ be a field of characteristic $2$ and $G=\mu_2$.  We define an action of $G$ on $V=kx\oplus ky$ as follows: for every $k$-scheme $T$ 
and every section $\zeta\in G(T)$, let $\zeta$ act on $V\otimes_k \OO_T$ by sending $x$ to $\zeta x$ and $y$ to $\zeta y$.  Then 
$\pi:\Vv\rightarrow\Vv/G$ is a $G$-torsor away from the one singular point in $\Vv/G$.  Hence, $\pi$ is 
ramified at every height 1 prime, but $G$ has no pseudo-reflections.
}
\end{example}
We must therefore take a different approach to showing that the action of $G/H$ on $\Vv/H$ has no pseudo-reflections at the origin.  Our strategy is 
to reduce to the classical case by lifting to characteristic 0.  This is carried out after some preliminary lemmas.
\begin{lemma}
\label{l:faithful}
Let $G$ be a finite group scheme which acts faithfully on an affine scheme $U$.  If $H$ is a normal subgroup scheme 
of $G$, then the action of $G/H$ on $U/H$ is faithful.
\end{lemma}
\begin{proof}
Let $\mf{X}=[U/H]$ and let $\pi:U\rightarrow U/H$ be the natural map.  We must show that if $G'$ is a subgroup scheme 
of $G$ such that $G'/H$ acts trivially on $U/H$, then $G'=H$.  Replacing $G$ by $G'$, we can assume $G'=G$.\\
\\
Since $G$ acts faithfully on $U$, there is a non-empty 
open substack of $\mf{X}$ which is isomorphic to its coarse space.  That is, we have a non-empty open subscheme $V$ of 
$U/H$ over which $\pi$ is an $H$-torsor.  Let $P=V\times_{U/H}U$.  Since $G$ acts on $P$ over $V$, we obtain a morphism 
\[
s:G\longrightarrow \mathcal{A}ut(P)=H.
\]
Note that $s$ is a section of the closed immersion $H\rightarrow G$, so $H=G$.
\end{proof}
\begin{lemma}
\label{l:polydvr}
Let $G$ be a finite flat linearly reductive group scheme over a complete discrete valuation ring $R$ with residue field $k$.  If $G$ acts linearly on $\mathbb{A}^n_R$ 
and $\mathbb{A}^n_k/G_k$ is isomorphic to $\mathbb{A}^n_k$, then $\mathbb{A}^n_R/G$ is isomorphic to $\mathbb{A}^n_R$.
\end{lemma}
\begin{proof}
Let $\mf{m}$ be the maximal ideal of $R$ and let $\mathbb{A}^n_R/G=\Spec A$.  
Since $\mathbb{A}^n_R$ is flat over $R$, it follows that $\mathbb{A}^n_R/G$ is as well (see \emph{e.g.} \cite[Thm 4.16(ix)]{alper}).  
Since $G$ is linearly reductive,
\[
\Spec k\times_R\mathbb{A}^n_R/G=\mathbb{A}^n_k/G_k.
\]
Choose an isomorphism 
\[
\varphi_0:k[x_1,\dots,x_n]\longrightarrow A\otimes_Rk
\]
and let $r_i\in R$ be an arbitrary lift of $\varphi_0(x_i)$.  By Nakayama's Lemma, the morphism 
\[
\varphi:R[x_1,\dots,x_n]\longrightarrow A
\]
sending $x_i$ to $r_i$ 
is surjective.  As $R$ is complete, to show $\varphi$ is an isomorphism, we need only show that the base change $\varphi_m$ of $\varphi$ to 
$R/\mf{m}^{\ell+1}$ is an isomorphism for every $\ell$.  This follows from the fact that $\varphi_0$ is an isomorphism and $A\otimes_R R/\mf{m}^\ell$ is flat over $R/\mf{m}^{\ell}$.
\end{proof}
\begin{proposition}
\label{prop:quotnopr}
Let $G$ be a finite linearly reductive group scheme over $S$ with a faithful action on a finite-dimensional $k$-vector 
space $V$.  Let $U=\Vv$ and $H$ be the subgroup scheme of $G$ generated by pseudo-reflections.  Then the induced 
action of $G/H$ on $U/H\simeq\mathbb{A}_k^n$ has no pseudo-reflections at the origin.
\end{proposition}
\begin{proof}
By the ``if'' direction of Theorem \ref{thm:newcst}, we have $k[V]^H=k[W]$ for some subvector space $W$ of $k[V]$.  
The proof of \cite[Prop 6.19]{neusel} shows that the degrees of the homogeneous generators of $k[V]^H$ are determined.  
As a result, the action of $G/H$ on $k[W]$ is linear.  Lemma \ref{l:faithful} further tells us that this action is 
faithful.\\
\\
Assume that the subgroup scheme $H''$ of $G/H$ generated by pseudo-reflections is non-trivial.  Then $H''=H'/H$ where $H'$ is a normal subgroup 
scheme of $G$ which properly contains $H$.  To prove $G/H$ has no pseudo-reflections at the origin, it suffices by Lemma \ref{l:prdescend} to 
replace $k$ by its algebraic closure.  
By \cite[Lemma 2.11]{tame}, we see then that $G$ is the semi-direct product of its identity component which is diagonalizable, 
and a finite constant tame group scheme.  The same is true for $H$ and $H'$.\\
\\
Let $R$ be a complete discrete valuation ring whose residue field is $k$ and whose fraction field $K$ is of characteristic 0.  
Lemma \ref{l:wsoverZ} shows that there exist 
finite flat linearly reductive group schemes $\til{G}$, $\til{H}$, and $\til{H}'$ over $R$ 
whose base change to $k$ are $G$, $H$, and $H'$, respectively.  Furthermore, $\til{H}'$ and $\til{H}$ are normal closed subgroup schemes 
of $\til{G}$, and $\til{H}$ is a proper subgroup scheme of $\til{H}'$.  
In characteristic 0, every finite flat group scheme is locally constant, so after replacing $R$ by a finite extension, we can further assume that 
$\til{G}_K$, $\til{H}_K$, and $\til{H}'_K$ are constant group schemes.\\
\\
Let $\mf{m}$ denote the maximal ideal of $R$ and let $R_\ell=R/\mf{m}^{\ell}$.  Let 
$\til{G}_{\ell}$,  $\til{H}_{\ell}$, and $\til{H}'_{\ell}$ denote the base change of $\til{G}$,  $\til{H}$, and $\til{H}'$ to $R_\ell$.  
Choosing a basis for $V$, we can identify $U$ with $\mathbb{A}^n_k$.  The $G$-action on $U$ is then given by a group scheme homomorphism 
$\varphi_0:G\longrightarrow GL_{n,k}$.  By \cite[Exp. III 2.3]{sga3}, given a deformation $\varphi_\ell:\til{G}_\ell\longrightarrow GL_{n,R_\ell}$ of 
$\varphi_0$, the obstruction to deforming $\varphi_\ell$ to a homomorphism $\varphi_{\ell+1}:\til{G}_{\ell+1}\longrightarrow GL_{n,R_{\ell+1}}$ lies in 
\[
H^2(\til{G}_{\ell},Lie(GL_n)\otimes\mf{m}^{\ell}/\mf{m}^{\ell+1}),
\]
which vanishes as $\til{G}_\ell$ is linearly reductive.  We therefore obtain a faithful action of $\til{G}$ on $\mathbb{A}^n_R$ lifting the action of 
$G$ on $U$.\\
\\
By Lemma \ref{l:polydvr}, we see that $\mathbb{A}^n_K/\til{H}_K$ and $\mathbb{A}^n_K/\til{H}'_K$ are polynomial.  
The classical Chevalley-Shephard-Todd theorem then shows that there is a pseudo-reflection 
$\til{N}_K$ of $\til{G}_K$ which is contained in $\til{H}'_K$ but not contained in $\til{H}_K$.  Note that this is not yet a contradiction as it is 
not clear that $\til{H}_K$ is the subgroup scheme of $\til{G}_K$ generated by pseudo-reflections.  
Let $\til{N}$ be the closure of $\til{N}_K$ in $\til{G}$.  
Since $\til{G}$ is a finite flat linearly reductive group scheme over $R$, we see that $\til{N}$ is as well.  Since $\til{N}_K$ is a pseudo-reflection, 
there exists some $v=\sum_ia_ix_i\in K[x_1,\dots,x_n]$ such $\til{N}_K$ acts trivially on $K[x_1,\dots,x_n]/v$.  After scaling the $a_i$, we can 
assume $a_1\in R^*$ and all $a_i\in R$.  Consider the commutative diagram
\[
\xymatrix{
0\ar[r] & vK[x_1,\dots,x_n]\ar[r] & K[x_1,\dots,x_n]\ar[r] & K[x_1,\dots,x_n]/v\ar[r] & 0\\
0\ar[r] & vR[x_1,\dots,x_n]\ar[r]\ar[u] & R[x_1,\dots,x_n]\ar[r]\ar[u] & R[x_1,\dots,x_n]/v\ar[r]\ar[u]^\psi & 0
}
\]
of $\til{N}$-comodules.  Since the left square is cartesian, we see that $\psi$ is injective.  It follows that the action of $\til{N}$ on the hyperplane 
defined by $v$ in $\mathbb{A}^n_R$ is trivial.  Reducing mod $\mf{m}$, we see that $\til{N}_k$ is a pseudo-reflection of $G$.  Furthermore, 
$\til{N}_k$ is not contained in $H$, which is a contradiction.
\end{proof}
Using Lemma \ref{l:faithful} and Proposition \ref{prop:quotnopr}, we prove the ``only if'' direction of Theorem \ref{thm:newcst}, assuming the special case 
of Theorem \ref{thm:torsor} in which $U=\Vv$.
\begin{proof}[Proof of ``only if'' direction of Theorem \ref{thm:newcst}]
Let $H$ be the subgroup scheme 
generated by pseudo-reflections.  By the ``if'' direction, $k[V]^H$ is polynomial and as explained in the proof of 
Proposition \ref{prop:quotnopr}, the $G/H$-action on $k[V]^H$ is linear.  Since $G/H$ acts faithfully on $U/H$ without 
pseudo-reflections at the origin by Lemma \ref{l:faithful} and Proposition \ref{prop:quotnopr}, and since $M=U/G$ is smooth by assumption, 
Theorem \ref{thm:torsor} implies that $U/H$ is a $G/H$-torsor over $U/G$ after potentially shrinking $U/G$.  Since the 
origin of $U/H$ is a fixed point, we conclude that $G=H$.
\end{proof}

\section{Theorem \ref{thm:torsor} for Linear Actions on Polynomial Rings}
\label{sec:polytorsor}
In Section \ref{sec:polynomial}, we reduced the proof of the ``only if'' direction of Theorem \ref{thm:newcst} to 
\begin{proposition}
\label{prop:polytorsor}
Let $G$ be a stable group scheme over $S$ which acts faithfully on a finite-dimensional 
$k$-vector space $V$.  Then Theorem \ref{thm:torsor} holds when $U=\Vv$ and $x$ is the origin.
\end{proposition}
The proof of this proposition is given in two steps.  We handle the case when $G$ is diagonalizable in 
Subsection \ref{subsec:iwanari} and then handle the general case in Subsection \ref{subsec:ptfinish} by making use of 
the diagonalizable case.

\subsection{Reinterpreting a Result of Iwanari}
\label{subsec:iwanari}
The key to proving Proposition \ref{prop:polytorsor} for diagonalizable $G$ is provided by Theorem 3.3 and Proposition 3.4 of 
\cite{mfr} after we reinterpret them in the language of pseudo-reflections.  We refer the reader to \cite[p.4-6]{mfr} for 
the basic definitions concerning monoids.  We recall the following definition given in \cite[Def 2.5]{mfr}.
\begin{definition}
\emph{An injective morphism $i:P\rightarrow F$ from a simplicially toric sharp monoid to a free monoid is called a 
\emph{minimal\ free\ resolution} if $i$ is close and if for all injective close morphisms $i':P\rightarrow F'$ to a 
free monoid $F'$ of the same rank as $F$, there 
is a unique morphism $j:F\rightarrow F'$ such that $i'=ji$.}
\end{definition}
Given a faithful action of a finite diagonalizable group scheme $\Delta$ over $S$ on a $k$-vector space $V$ of dimension $n$, we can 
decompose $V$ as a direct sum of 1-dimensional $\Delta$-representations.  Therefore, after choosing an appropriate basis, 
we have an identification of $k[V]$ with $k[\N^n]$ and can assume that the $\Delta$-action on $U=\Vv$ is induced from a 
morphism of monoids
\[
\pi:F=\N^n\longrightarrow A,
\]
where $A$ is the finite abelian group such that $\Delta$ is the Cartier dual $D(A)$ of $A$.  We see then that 
\[
U/\Delta=\Spec k[P],
\]
where $P$ is the submonoid $\{p\mid\pi(p)=0\}$ of $F$.  Note that $P$ is simplicially toric sharp, that 
$i:P\rightarrow F$ is close, and that $A=F^{gp}/i(P^{gp})$.\\
\\
We now give the relationship between minimal free resolutions and pseudo-reflections.
\begin{proposition}
\label{prop:mfrpr}
With notation as above, 
$i:P\rightarrow F$ is a minimal free resolution if and only if the action of $\Delta$ on $V$ has no pseudo-reflections.
\end{proposition}
\begin{proof}
If $i$ is not a minimal free resolution, then without loss of generality, $i=ji'$, where $i':P\rightarrow F$ is close and 
injective, and $j:F\rightarrow F$ is given by
\[
j(a_1,a_2,\dots, a_n)=(ma_1,a_2,\dots, a_n)
\]
with $m\neq 1$.  We have then a short exact sequence 
\[
0\longrightarrow F^{gp}/i'(P^{gp})\longrightarrow F^{gp}/i(P^{gp})\longrightarrow F^{gp}/(m,1,\dots,1)(F^{gp})
\longrightarrow 0.
\]
Let $N$ be the Cartier dual of $F^{gp}/(m,1,\dots,1)(F^{gp})$, which is a subgroup scheme of $\Delta$.  Letting 
$\{x_i\}$ be the standard basis of $F$, we see that
\[
k[F]^N=k[x_1^m,x_2,\dots,x_n],
\]
and so $V^N$, which is the 
degree 1 part of $k[F]^N$, has codimension 1 in $V$.  Therefore, $N$ is a pseudo-reflection.\\
\\
Conversely, suppose $N$ is a pseudo-reflection.  Since $N$ is a subgroup scheme of $\Delta$, it is diagonalizable 
as well.  Let $N=\Spec k[B]$, where $B$ is a finite abelian group and let $\psi:A\rightarrow B$ be the induced map.  
We see that
\[
V^N=\bigoplus_{i\neq j}kx_i
\] 
for some $j$.  Without loss of generality, $j=1$.  It follows then that
\[
\{f\in F\mid\psi\pi(f)=0\}=(m,1,\dots,1)F
\]
for some $m$ dividing $|B|$.  Since the $\Delta$ action on $V$ is assumed to be faithful, we 
see, in fact, that $m=|B|$.  Therefore, $i$ factors through $\cdot(m,1,\dots,1):F\longrightarrow F$, which 
shows that $i$ is not a minimal free resolution.
\end{proof}
Having reinterpreted minimal free resolutions, the proof of Proposition \ref{prop:polytorsor} for 
diagonalizable group schemes $G$ follows easily from Iwanari's work.
\begin{proposition}
\label{prop:ptdiag}
Let $G=\Delta$ be a finite diagonalizable group scheme over $S$ which acts faithfully on a finite-dimensional $k$-vector space $V$.  
Then Theorem \ref{thm:torsor} holds when $U=\Vv$ and $x$ is the origin.  In 
this case it is not necessary to shrink $M$ to a smaller Zariski neighborhood of 
the image of $x$.
\end{proposition}
\begin{proof}
Let $F$ and $P$ be as above, and let $\mf{X}=[U/\Delta]$.  By 
Proposition \ref{prop:mfrpr}, 
the morphism $i:P\rightarrow F$ is a minimal free resolution.  Theorem 3.3 (1) and Proposition 3.4 of \cite{mfr} then 
show that the natural morphism $\mf{X}\times_M M^0\rightarrow M^0$ is an isomorphism.  Since $\mf{X}\times_M M^0=
[U^0/\Delta]$, we see $U^0$ is a $\Delta$-torsor over $M^0$.
\end{proof}

\subsection{Finishing the Proof}
\label{subsec:ptfinish}
The goal of this subsection is to prove Proposition \ref{prop:polytorsor}.  
The main result used in the proof of this proposition, as well as in the proof of Theorem \ref{thm:torsor}, is the following.
\begin{proposition}
\label{prop:et}
Let notation and hypotheses be as in Theorem \ref{thm:torsor}.  Let $X=U/\Delta$ and $G=\Delta\rtimes Q$, where  
$\Delta$ is diagonalizable and $Q$ is constant and tame.  If in addition to assuming that 
$G$ acts without pseudo-reflections at $x$, we assume that $\Delta$ is local and that the base change of $U$ to 
$X^{sm}$ is a $\Delta$-torsor over $X^{sm}$, then after possibly shrinking $M$ to a smaller Zariski neighborhood of 
the image of $x$, 
the quotient map $f:X\rightarrow M$ is unramified in codimension 1.
\end{proposition}
\begin{proof}
Let $g$ be the quotient map $U\rightarrow X$.  
For every $q\in Q$, consider the cartesian diagram
\[
\xymatrix{
Z_q\ar[r]\ar[d] & U\ar[d]^{\Delta}\\
U\ar[r]^-{\Gamma_q} & U\times U
}
\]
where $\Gamma_q(u)=(u,qu)$.  We see that $Z_q$ is a closed subscheme of $U$ and that $Z_q(T)$ is the set of 
$u\in U(T)$ which are fixed by $q$.  Let $Z$ be the closed subset of $U$ which is the union of the $Z_q$ for $q\neq 1$.  
Since the action of $G$ on $U$ is faithful, $Z$ is not all of $U$.  Let $Z'$ be the union of the 
codimension 1 components of $Z$.  Since $fg$ is finite, we see that $fg(Z')$ is a closed subset of $M$.  Moreover, 
$fg(Z')$ does not contain 
the image of $x$, as $G$ is assumed to act without pseudo-reflections at $x$.  By shrinking 
$M$ to $M-fg(Z')$, we can assume that no non-trivial $q\in Q$ acts trivially on a divisor of $U$.\\
\\
Let $U=\Spec R$.  
The morphism $f$ is unramified in codimension 1 if and only if the (traditional) inertia groups of all height 1 primes 
$\mf{p}$ of $R^{\Delta}$ are trivial.  
So, we must show that if $q\in Q$ acts trivially on $V(\mf{p})$, then $q=1$.  Since $g$ is finite, and 
hence integral, the going up theorem shows that 
\[
\mf{p}R=\mf{P}_1^{e_1}+\dots +\mf{P}_n^{e_n},
\]
where the $\mf{P}_i$ are height 1 primes and the $e_i$ are positive integers.  
Note that $X$ is normal and so the complement of 
$X^{sm}$ in $X$ has codimension at least 2.  As a result, 
\[
h:U\times_X \Spec \OO_{X,\mf{p}}\longrightarrow \Spec \OO_{X,\mf{p}}
\]
is a $\Delta$-torsor.  Since $\Delta$ is local, $h$ is a homeomorphism of topological spaces, so 
there is exactly one prime $\mf{P}$ lying over $\mf{p}$.  We see then that $U\times_X V(\mf{p})=V(\mf{P}^e)$ for 
some $e$.\\
\\
Let $V(\mf{p})^0$ be the intersection of $V(\mf{p})$ with $X^{sm}$, and let $Z^0=U\times_X V(\mf{p})^0$.  
Then $Z^0$ is a $\Delta$-torsor over $V(\mf{p})^0$.  Since $q$ acts trivially on $V(\mf{p})$, we obtain an action of 
$q$ on $Z^0$ over $V(\mf{p})^0$, and hence a group scheme homomorphism
\[
\varphi:Q'_{V(\mf{p})^0}\longrightarrow \mathcal{A}ut(Z^0/V(\mf{p})^0)=\Delta_{V(\mf{p})^0},
\]
where $Q'$ denotes the subgroup of $Q$ generated by $q$.  Since $V(\mf{p})^0$ is reduced, we see that $\varphi$ 
factors through the reduction of $\Delta_{V(\mf{p})^0}$, which is the trivial group scheme.  Therefore, $q$ acts 
trivially on $Z^0$.\\
\\
Since the complement of $X^{sm}$ in $X$ has codimension at least 2, and since $g$ factors as a 
flat map $U\rightarrow [U/\Delta]$ followed by a coarse space map $[U/\Delta]\rightarrow X$, 
both of which are codimesion-preserving (see Definition 4.2 and Remark 4.3 of \cite{fant}), we see that the complement of 
$Z^0$ in $V(\mf{P}^e)$ has codimension at least 2.  
Note that if $Y$ is a normal scheme and $W$ is an open subscheme of $Y$ whose complement has 
codimension at least 2, then any morphism from $W$ to an affine scheme $Z$ extends uniquely to a morphism from $Y$ to 
$Z$.  Since the action of $q$ on $V(\mf{P}^e)$ restricts to a trivial action on $Z^0$, the action of 
$q$ on $V(\mf{P}^e)$ is trivial.  Therefore, $q$ acts trivially on a divisor of $U$, and so $q=1$.
\end{proof}
\begin{proof}[Proof of Proposition \ref{prop:polytorsor}]
Let $k'/k$ be a finite Galois extension such that $G_{k'}\simeq\Delta\rtimes Q$, where $\Delta$ 
is diagonalizable and $Q$ is constant and tame.  Let $S'=\Spec k'$ and consider the diagram 
\[
\xymatrix{
U'\ar[r]\ar[d] & U\ar[d]\\
M'\ar[r]\ar[d] & M\ar[d]\\
S'\ar[r] & S
}
\]
where the squares are cartesian.  We denote by $x'$ the induced $k'$-rational point of $U'$.  
Since $\Delta$ is the product of a local diagonalizable group scheme 
and a locally constant diagonalizable group scheme, replacing $k'$ by a further extension if necessary, we can assume 
that $\Delta$ is local.\\
\\
Since $G$ is stable, $G_{k'}$ has no pseudo-reflections at $x'$.  It follows then from Proposition 
\ref{prop:et} that there exists an open neighborhood $W'$ of $x'$ such that $U'\times_{M'}W'\longrightarrow W'$ 
is unramified in codimension 1.  Since $k'/k$ is a finite Galois extension, replacing $W'$ by the intersection of the 
$\tau(W')$ as $\tau$ ranges over the elements of $\Gal(k'/k)$, we can assume $W'$ is Galois-invariant.  Hence,
$W'=W\times_M M'$ for some open subset $W$ of $M$.  We shrink $M$ to $W$.\\
\\
To check that $U^0$ is a $G$-torsor over $M^0$, we can look \'etale locally.  We can therefore assume $S=S'$.  
Let $X=U/\Delta$, and let $g:U\rightarrow X$ and $f:X\rightarrow M$ be the quotient maps.  We denote by $X^0$ the 
fiber product $X\times_M M^0$ and by $f^0$ the induced morphism $X^0\rightarrow M^0$.\\
\\
By Proposition \ref{prop:ptdiag}, we know that the base change of $U$ to $X^{sm}$ is a $\Delta$-torsor over $X^{sm}$.  
Since $f$ is unramified in codimension 1, we see that $f^0$ is as well.  Since $M^0$ is smooth and 
$X^0$ is normal, the purity of the branch locus theorem \cite[X.3.1]{sga1} implies that $f^0$ is \'etale, 
and hence a $Q$-torsor.  Since $X^0$ is \'etale over $M^0$, it is smooth.  As a result, $U^0$ is a $\Delta$-torsor over 
$X^0$ from which it follows that $U^0$ is a $G$-torsor over $M^0$.
\end{proof}
This finishes the proof of Proposition \ref{prop:polytorsor}, and hence also of Theorem \ref{thm:newcst}.  We conclude 
this section by proving Corollary \ref{cor:newcst}.
\begin{proof}[Proof of Corollary \ref{cor:newcst}]
Let $U=\Spec R$ and $M=U/G$.  We denote by $y$ the image of $x$.  Since $G$ being generated by pseudo-reflections at $x$ 
implies that $G_K$ is generated by pseudo-reflections at $x$ for arbitrary finite linearly reductive group schemes $G$, and since 
smoothness of $M$ at $y$ can be checked \'etale locally, we can assume that $x$ is $k$-rational.  
Let $V=\mf{m}_x/\mf{m}_x^2$ be the cotangent space of $x$.  As $G$ is linearly reductive, there is a $G$-equivariant section of 
$\mf{m}_x \rightarrow V$.  This yields a $G$-equivariant map $\Sym^\bullet(V)\rightarrow R$, which induces an isomorphism 
$k[[V]]\longrightarrow \hat{\OO}_{U,x}$ of $G$-representations.  That is, complete locally, we have 
linearized the $G$-action.  Since $\hat{\OO}_{M,y}=k[[V]]^G$, the corollary follows from Theorem \ref{thm:newcst},
as $M$ is smooth at $y$ if and only if $\hat{\OO}_{M,y}$ is a formal power series ring over $k$.
\end{proof}

\section{Actions on Smooth Schemes}
\label{sec:torsor}
Having proved Theorem \ref{thm:torsor} for polynomial rings with linear actions, we now turn to the general case.  
We begin with two preliminary lemmas and a technical proposition.
\begin{lemma}
\label{l:repable}
Let $U$ be a smooth affine scheme over $S$ with an action of a finite diagonalizable group scheme $\Delta$.  
Then there is a closed subscheme $Z$ of $U$ on which $\Delta$ acts trivially, and with the property that  
every closed subscheme $Y$ on which $\Delta$ acts trivially factors through $Z$.  Furthermore, the construction of 
$Z$ commutes with flat base change on $U/\Delta$.
\end{lemma}
\begin{proof}
Let $U=\Spec R$ and $\Delta=\Spec k[A]$, where $A$ is a finite abelian group written additively.  The $\Delta$-action on 
$U$ yields an $A$-grading
\[
R=\bigoplus_{a\in A}R_a.
\]
We see that if $\mathcal{J}$ is an ideal of $R$, then $\Delta$ 
acts trivially on $Y=\Spec R/\mathcal{J}$ if and only if $\mathcal{J}$ contains the $R_a$ for $a\neq 0$.  Letting 
$\mathcal{I}$ be the ideal generated by the $R_a$ for $a\neq0$, we see that $\Spec R/\mathcal{I}$ is our desired $Z$.\\
\\
We now show that the formation of $Z$ commutes with flat base change.  Note that
\[
U/\Delta=\Spec R_0.
\]
Let $R'_0$ be a 
flat $R_0$-algebra and let $R'=R'_0\otimes_{R_0}R$.  The induced $\Delta$-action on $\Spec R'$ corresponds to the 
$A$-grading
\[
R'=\bigoplus_{a\in A}(R'_0\otimes_{R_0}R_a).
\]
Since $R'_0$ is flat over $R_0$, we see that 
$\mathcal{I}\otimes_{R_0}R'_0$ is an ideal of $R'$, and one easily shows that it is the ideal generated by the 
$R'_0\otimes_{R_0}R_a$ for $a\neq0$.
\end{proof}
Recall that if $G$ is a group scheme over a base scheme $B$ which acts on a $B$-scheme $U$, and if $y:T\rightarrow U$ 
is a morphism of $B$-schemes, then the stabilizer group scheme $G_y$ is defined by the cartesian diagram
\[
\xymatrix{
G_y\ar[r]\ar[d] & G\times_B U\ar[d]^\varphi\\
T\ar[r]^-{y\times y} & U\times_B U
}
\]
where $\varphi(g,u)=(gu,u)$.  If $U$ is separated over $B$, then $G_y$ is a closed subgroup scheme of $G_T$.
\begin{lemma}
\label{l:fieldtorsor}
Let $B$ be a scheme and $G$ a finite flat group scheme over $B$.  If $G$ acts on a $B$-scheme $U$, then 
$U\rightarrow U/G$ is a $G$-torsor if and only if the stabilizer group schemes $G_y$ are trivial for all closed 
points $y$ of $U$.
\end{lemma}
\begin{proof}
The ``only if'' direction is clear.  To prove the ``if'' direction, it suffices to show that the stabilizer group 
schemes $G_y$ are trivial 
for all scheme valued points $y:T\rightarrow U$.  This is equivalent to showing that the universal stabilizer $G_u$ is 
trivial, where $u:U\rightarrow U$ is the identity map.  Since $G_u$ is a finite group scheme over $U$, it is given by 
a coherent sheaf $\mathcal{F}$ on $U$.  The support of $\mathcal{F}$ is a closed subset, and so to prove 
$G_u$ is trivial, it suffices to check this on stalks of closed points.  Nakayama's Lemma then shows that we need only 
check the triviality of $G_u$ on closed fibers.  That is, we need only check that the $G_y$ are trivial for closed 
points $y$ of $U$.
\end{proof}
\begin{proposition}
\label{prop:artin}
Let $U$ be a smooth affine scheme over $S$ with a faithful action of a stable group scheme $G$ fixing 
a $k$-rational point $x$.  If $N$ has a pseudo-reflection at $x$, then there is an \'etale neighborhood $T\longrightarrow U/G$ of $x$ and a divisor $D$ of $U_T$ defined 
by a principal ideal on which $N_T$ acts trivially.
\end{proposition}
\begin{proof}
Let $M=U/G$ and let $y$ be the image of $x$ in $M$.  As in the proof of Corollary \ref{cor:newcst}, we have an 
isomorphism $k[[V]]\longrightarrow \hat{\OO}_{U,x}$ of $G$-representations, where $V=\mf{m}_x/\mf{m}_x^2$.  If $N$ 
is a pseudo-reflection at $x$, then there is some $v\in V$ such that $N$ acts trivially on the closed subscheme of 
$\Spec k[[V]]$ defined by the prime ideal generated by $v$.\\
\\
Consider the contravariant functor $F$ which sends an $M$-scheme $T$ to the set of 
divisors of $U_T$ defined by a principal ideal on which $N_T$ acts trivially.  
As $F$ is locally of finite presentation and 
$U\times_M \Spec\hat{\OO}_{M,y}=\Spec\hat{\OO}_{U,x}$, Artin's Approximation Theorem \cite{approx} finishes the proof.
\end{proof}
We are now ready to prove Theorem \ref{thm:torsor}.  Our method of proof is similar to that of 
Proposition \ref{prop:polytorsor}; we first prove the theorem in the case that $G$ is diagonalizable and then make 
use of this case to prove the theorem in general.
\begin{proposition}
\label{prop:smt}
Theorem \ref{thm:torsor} holds when $G=\Delta$ is a finite diagonalizable group scheme.
\end{proposition}
\begin{proof}
Let $g:U\rightarrow M$ be the quotient map.  
Since any subgroup 
scheme $N$ of $\Delta$ is again finite diagonalizable, Lemma \ref{l:repable} shows that for every $N$, there exists a 
closed subscheme $Z_N$ of $U$ on which $N$ acts trivially, and with the 
property that every closed subscheme $Y$ on which $N$ acts trivially factors through $Z_N$.  Let $Z$ be the union of 
the finitely many closed subsets $Z_N$ for $N\neq 1$.  Since the action of 
$\Delta$ on $U$ is faithful, $Z$ has codimension at least 1.  Let $Z'$ be the union of all irreducible components of $Z$ 
which have codimension 1.  Since $\Delta$ acts without pseudo-reflections at $x$, we see $x\notin Z'$.  Note that $g(Z')$ 
is closed as $g$ is proper.  Since the construction of $Z$ commutes with flat base change on $M$ and since flat morphisms 
are codimension-preserving, replacing $M$ with $M-g(Z')$, we can assume that  
there are no non-trivial subgroup schemes of $\Delta$ which fppf locally on 
$M$ act trivially on a divisor of $U$.\\
\\
By Lemma \ref{l:fieldtorsor}, to show $U^0$ is a $\Delta$-torsor over $M^0$, it suffices to show that for every 
closed point $y$ of $U$ which maps to $M^0$, the stabilizer group scheme $\Delta_y$ is trivial.  Fix such a closed 
point $y$ and let $T=\Spec k(y)$.  Since $T$ is fppf over $S$, we see from Proposition \ref{prop:artin} that the closed 
subgroup scheme $\Delta_y$ of $\Delta_T$ acts faithfully on $U_T$ without pseudo-reflections 
at the $k(y)$-rational point $y'$ of $U_T$ induced by $y$.  Since $y$ maps to a smooth point of $M$, it follows that 
$y'$ maps to a smooth point of $M_T$.  Corollary \ref{cor:newcst} then shows that $\Delta_y$ is generated by 
pseudo-reflections.  Since $\Delta_y$ has no pseudo-reflections, it is therefore trivial.
\end{proof}
\begin{proof}[Proof of Theorem \ref{thm:torsor}]
If $G=\Delta\rtimes Q$, where $\Delta$ is diagonalizable and $Q$ is constant and tame, then letting $Z'$ be as in 
Proposition \ref{prop:smt} and letting $U$, $X$, $f$, and $g$ be as in the proof of Proposition \ref{prop:polytorsor}, 
the proof of Proposition \ref{prop:smt} shows that after replacing $M$ by $M-fg(Z')$, the base change of $U$ to $X^{sm}$ 
is a $\Delta$-torsor over $X^{sm}$.  As in the proof of Proposition \ref{prop:polytorsor}, we can then reduce the 
general case to the case when $G=\Delta\rtimes Q$, where $\Delta$ is local diagonalizable and $Q$ is constant tame.  
The last paragraph of the proof of Proposition \ref{prop:polytorsor} then shows that $U^0$ is a $G$-torsor over $M^0$.
\end{proof}

\section{Schemes with Linearly Reductive Singularities}
\label{sec:lrs}
Let $k$ be a perfect field of characteristic $p$.
\begin{definition}
\emph{We say a scheme $M$ over $S$ has $\emph{linearly}$ $\emph{reductive}$ $\emph{singularities}$ if 
there is an \'etale cover $\{U_i/G_i\rightarrow M\}$, where the $U_i$ are smooth over $S$ and the 
$G_i$ are linearly reductive group schemes which are finite over $S$.}
\end{definition}
Note that if $M$ has linearly reductive singularities, then it is automatically normal and in fact Cohen-Macaulay by 
\cite[p.115]{cm}.\\
\\
Our goal in this section is to prove Theorem \ref{prop:cslrs}, 
which generalizes the result that every scheme with quotient singularities prime to the characteristic is the coarse 
space of a smooth Deligne-Mumford stack.  
We remark that in the case of quotient singularities, the converse of the 
analogous theorem is true as well; that is, every scheme which is the coarse 
space of a smooth Deligne-Mumford stack has quotient singularities.  It is not clear, however, that the converse of 
Theorem \ref{prop:cslrs} should hold.  We know from Theorem 3.2 of \cite{tame} that $\mf{X}$ is \'etale locally 
$[V/G_0]$, where $G_0$ is a finite flat linearly reductive group scheme over $V/G_0$, but $V$ need not be smooth and 
$G_0$ need not be the base change of a group scheme over $S$.  On the other hand, Proposition \ref{l:semidirect} below 
shows that $\mf{X}$ is \'etale locally $[U/G]$ where $U$ is smooth and $G$ is a group scheme over $S$, but here $G$ is 
not finite.\\
\\
Before proving Theorem \ref{prop:cslrs}, we begin with a technical proposition 
followed by a series of lemmas.
\begin{proposition}
\label{l:semidirect}
Let $\mf{X}$ be a tame stack over $S$ with coarse space $M$.  Then there exists an \'etale cover $T\rightarrow M$ such that 
\[
\mf{X}\times_M T = [U/\mathbb{G}_{m,T}^r\rtimes H],
\]
where $H$ is a finite constant tame group scheme and $U$ is affine over $T$.  Furthermore, $\mathbb{G}_{m,T}^r\rtimes H$ 
is the base change to $T$ of a group scheme $\mathbb{G}_{m,S}^r\rtimes H$ over $S$, so $\mf{X}\times_M T= 
[U/\mathbb{G}_{m,S}^r\rtimes H]$.
\end{proposition}
\begin{proof}
Theorem 3.2 of \cite{tame} shows that there exists an \'etale cover $T\rightarrow M$ and a finite flat linearly reductive 
group scheme $G_0$ over $T$ acting on a finite finitely presented scheme $V$ over $T$ such that
\[
\mf{X}\times_M T=[V/G_0].
\]
By \cite[Lemma 2.20]{tame}, after replacing $T$ by a finer \'etale cover if necessary, we can assume there is a short exact 
sequence
\[
1 \longrightarrow\Delta \longrightarrow G_0 \longrightarrow H \longrightarrow 1,
\]
where $\Delta=\Spec \OO_T[A]$ is a finite diagonalizable group scheme and $H$ is a finite constant tame group scheme.  Since 
$\Delta$ is abelian, the conjugation action of $G_0$ on $\Delta$ passes to an action
\[
H\longrightarrow \Aut(\Delta)=\Aut(A).
\]
Choosing a surjection $F\rightarrow A$ in the category of $\Z[H]$-modules from a free module $F$, yields an $H$-equivariant 
morphism $\Delta\hookrightarrow \mathbb{G}_{m,T}^r$.  Using the $H$-action on $\mathbb{G}_{m,T}^r$, we define the group 
scheme $\mathbb{G}_{m,T}^r\rtimes G_0$ over $T$.  Note that there is an embedding
\[
\Delta\hookrightarrow \mathbb{G}_{m,T}^r\rtimes G_0
\]
sending $\delta$ to $(\delta,\delta^{-1})$, which realizes $\Delta$ as a normal subgroup scheme of 
$\mathbb{G}_{m,T}^r\rtimes G_0$.  We can therefore define
\[
G:=(\mathbb{G}_{m,T}^r\rtimes G_0)/\Delta.
\]
One checks that there is a commutative diagram
\[
\xymatrix{
1\ar[r] & \Delta\ar[r]\ar[d] & G_0\ar[r]\ar[d] & H\ar[r]\ar[d]^{id} & 1\\
1\ar[r] & \mathbb{G}_{m,T}^r\ar[r] & G\ar[r]^{\pi} & H\ar[r] & 1
}
\]
where the rows are exact and the vertical arrows are injective.\\
\\
We show that \'etale locally on $T$, there is a group scheme-theoretic section of $\pi$, so that 
$G=\mathbb{G}_{m,T}^r\rtimes H$.  Let $P$ be the sheaf on $T$ such that for any $T$-scheme $W$, $P(W)$ is the set of group 
scheme-theoretic sections of $\pi_W:G_W\rightarrow H_W$.  Note that the sheaf $\underline{\Hom}(H,G)$ parameterizing group 
scheme homomorphisms from $H$ to $G$ is representable since it is a closed subscheme of $G^{\times |H|}$ cut out by suitable 
equations.  We see that $P$ is the equalizer of the two maps
\[
\xymatrix{
\underline{\Hom}(H,G) \ar@<+.5ex>[r]^-{p_1} \ar@<-.5ex>[r]_-{p_2} & H^{\times |H|}
}
\]
where $p_1(\phi)=(\pi\phi(h))_h$ and $p_2(\phi)=(h)_h$.  That is, there is a cartesian diagram
\[
\xymatrix{
P\ar[r]\ar[d] & \underline{\Hom}(H,G)\ar[d]^{(p_1,p_2)}\\
H^{\times |H|}\ar[r]^-{\Delta} & H^{\times |H|}\times H^{\times |H|}
}
\]
Since $H$ is separated over $T$, we see that $P$ is a closed subscheme of $\underline{\Hom}(H,G)$.  In particular, it is 
representable and locally of finite presentation over $T$.  Furthermore, $P\rightarrow T$ is surjective as 
\cite[Lemma 2.16]{tame} shows that it has a section fppf locally.  To show $P$ has a section \'etale 
locally, by \cite[17.16.3]{ega4}, it suffices to prove $P$ is smooth over $T$.\\
\\
Given a commutative diagram
\[
\xymatrix{
X_0=\Spec A/\mathcal{I}\ar[r]\ar[d] & P\ar[d]\\
X=\Spec A\ar[r]\ar@{-->}[ur] & T
}
\]
with $\mathcal{I}$ a square zero ideal, we want to find a dotted arrow making the diagram commute.  That is, given a 
group scheme-theoretic section $s_0:G_{W_0}\rightarrow H_{W_0}$ of $\pi_{W_0}$, we want to find a group scheme homomorphism 
$s:G_W\rightarrow H_W$ which pulls back to $s_0$ and such that $\pi_W\circ s$ is the identity.  Note first that any group 
scheme homomorphism $s$ which pulls back to $s_0$ is automatically a section of $\pi_W$ since $H$ is a finite constant group 
scheme and $\pi_W\circ s$ pulls back to the identity over $W_0$.  By \cite[Exp. III 2.3]{sga3}, the obstruction to lifting 
$s_0$ to a group scheme homomorphism lies in
\[
H^2(H,Lie(G)\otimes\mathcal{I}),
\]
which vanishes as $H$ is linearly reductive.  This proves the smoothness of $P$.\\
\\
To complete the proof of the lemma, let $U:=V\times^{G_0}G$ and note that
\[
\mf{X}\times_M T=[V/G_0]=[U/G].
\]
Since $V$ is finite over $T$ and $G$ is affine over $T$, it follows that $U$ is affine over $T$ as well.  Replacing $T$ by a 
finer \'etale cover if necessary, we have
\[
\mf{X}\times_M T=[U/\mathbb{G}_{m,T}^r\rtimes H].
\]
Lastly, 
the scheme underlying $\mathbb{G}_{m,T}^r\rtimes H$ is $\mathbb{G}_{m,T}^r\times_T H$ and its group scheme 
structure is determined by the action $H\rightarrow\Aut(\mathbb{G}_{m,T}^r)$.  Since $\Aut(\mathbb{G}_{m,T}^r)=\Aut(\Z^r)$, 
we can use this same action to define the semi-direct product $\mathbb{G}_{m,S}^r\rtimes H$ and it is clear that this group 
scheme base changes to $\mathbb{G}_{m,T}^r\rtimes H$.
\end{proof}
\begin{lemma}
\label{l:smoothstack}
If $V$ is a smooth $S$-scheme with an action of finite linearly reductive group scheme $G_0$ over $S$, then $[V/G_0]$ is 
smooth over $S$.
\end{lemma}
\begin{proof}
Let $\mf{X}=[V/G_0]$.  To prove $\mf{X}$ is smooth, it suffices to work \'etale locally on $S$, where, by 
\cite[Lemma 2.20]{tame}, we can assume $G_0$ fits into a short exact sequence 
\[
1 \longrightarrow\Delta \longrightarrow G_0 \longrightarrow H \longrightarrow 1,
\]
where $\Delta$ is a finite diagonalizable group scheme and $H$ is a finite constant tame group scheme.  Let $G$ be 
obtained from $G_0$ as in the proof of Proposition \ref{l:semidirect} and let $U=V\times^{G_0}G$.  Since $\mf{X}=[U/G]$, it 
suffices to show 
$U$ is smooth over $S$.  The action of $G_0$ on $V\times G$, given by $g_0\cdot (v,g)=(vg_0,g_0g)$, is free as the 
$G_0$-action on $G$ is free.  As a result, $U=[(V\times G)/G_0]$ and $G/G_0=[G/G_0]$.  Since the projection map 
$p:V\times G\rightarrow G$ is $G_0$-equivariant, we have a cartesian diagram
\[
\xymatrix{
V\times G\ar[r]^p \ar[d] & G\ar[d]\\
U\ar[r]^q & G/G_0
}
\]
Since $p$ is smooth, $q$ is as well.  Since $G\rightarrow [G/G_0]=G/G_0$ is flat and $G$ is smooth, \cite[17.7.7]{ega4} shows 
that $G/G_0$ is smooth, and so $U$ is as well.
\end{proof}
\begin{lemma}
\label{l:torsor}
Let $X$ be a smooth $S$-scheme and $i:U\hookrightarrow X$ an open subscheme whose complement has codimension at least 2.  
Let $P$ be a $G$-torsor on $U$, where $G=\mathbb{G}_m^r \rtimes H$ and $H$ is a finite constant \'etale group scheme.  
Then $P$ extends uniquely to a $G$-torsor on $X$.
\end{lemma}
\begin{proof}
The structure map from $P$ to $U$ factors as $P\rightarrow P_0\rightarrow U$, 
where $P$ is a $\mathbb{G}_m^r$-torsor over $P_0$ and $P_0$ is an $H$-torsor over $U$.  Since the complement of $U$ in $X$ 
has codimension at least $2$, we have 
$\pi_1(U)=\pi_1(X)$ and so $P_0$ extends uniquely to an $H$-torsor $Q_0$ on $X$.  Let $i_0:P_0\hookrightarrow Q_0$ be the 
inclusion map.  Since $Q_0$ is smooth and the complement of $P_0$ in $Q_0$ has codimension at least 2, the natural map 
$\textrm{Pic}(Q_0)\rightarrow\textrm{Pic}(P_0)$ is an isomorphism.  It follows that any line bundle over $P_0$ can be 
extended uniquely to a line bundle over $Q_0$.  We can therefore inductively construct a unique lift of $P$ over $X$.
\end{proof}
Our proof of the following lemma closely follows that of \cite[Thm 4.6]{fant}.
\begin{lemma}
\label{l:univprop}
Let $f:\mathcal{Y}\rightarrow M$ be an $S$-morphism from a smooth tame stack $\mathcal{Y}$ to its coarse space which 
pulls back to an isomorphism over the smooth locus $M^0$ of $M$.  If $h:\mf{X}\rightarrow M$ is a dominant, 
codimension-preserving morphism \emph{(}see \emph{\cite[Def 4.2]{fant})} from a smooth tame stack, then there is a morphism 
$g:\mf{X}\rightarrow\mathcal{Y}$, unique up to unique isomorphism, such that $fg=h$.
\end{lemma}
\begin{proof}
We show that if such a morphism $g$ exists, then it is unique.  Suppose $g_1$ and $g_2$ are two such morphisms.  We see 
then that $g_1|_{h^{-1}(M^0)}=g_2|_{h^{-1}(M^0)}$.  Since $h$ is dominant and codimension-preserving, $h^{-1}(M^0)$ is 
open and dense in $\mf{X}$.  Proposition 1.2 of \cite{fant} shows that if $\mf{X}$ and $\mathcal{Y}$ are 
Deligne-Mumford with $\mf{X}$ normal and $\mathcal{Y}$ separated, then $g_1$ and $g_2$ are uniquely isomorphic.  The proof, 
however, applies equally well to tame stacks since the only key ingredient used about Deligne-Mumford stacks is that they 
are locally $[U/G]$ where $G$ is a separated group scheme.\\
\\
By uniqueness, to show the existence of $g$, we can assume by Proposition \ref{l:semidirect} that 
$\mathcal{Y}=[U/G]$, where $U$ is smooth and affine, and 
$G=\mathbb{G}_m^r \rtimes H$, where $H$ is a finite constant tame group scheme.  Let $p:V\rightarrow\mf{X}$ be a smooth 
cover by a smooth scheme.  Since smooth morphisms are dominant and codimension-preserving, uniqueness implies that 
to show the existence of $g$, we need only show there is a morphism 
$g_1:V\rightarrow\mathcal{Y}$ such that $fg_1=hp$.  So, we can assume $\mf{X}=V$.\\
\\
Given a stack $\mathcal{Z}$ over $M$, let $\mathcal{Z}^0=M^0\times_M \mathcal{Z}$.  Given a morphism 
$\pi:\mathcal{Z}_1\rightarrow\mathcal{Z}_2$ of $M$-stacks, let $\pi^0:\mathcal{Z}_1^0\rightarrow\mathcal{Z}_2^0$ denote 
the induced morphism.  Since $f^0$ is an isomorphism, there is a morphism 
$g^0:V^0\rightarrow\mathcal{Y}^0$ such that $f^0g^0=h^0$.  It follows that there is a $G$-torsor $P^0$ over $V^0$ and a 
$G$-equivariant map from $P^0$ to $U^0$ such that the diagram 
\[
\xymatrix{
P^0\ar[r]\ar[d] & U^0\ar[d]\\
V^0\ar[d]\ar[r] & \mathcal{Y}^0\ar[dl]^{\simeq}\\
M^0
}
\]
commutes and the square is cartesian.  By Lemma \ref{l:torsor}, $P^0$ extends to a $G$-torsor $P$ over $V$.\\
\\
Note that if $X$ is a normal algebraic space and $i:W\hookrightarrow X$ is an open subalgebraic space whose complement has 
codimension at least 2, then any morphism from $W$ to an affine scheme $Y$ extends uniquely to a morphism 
$X\rightarrow Y$.  As a result, the morphism from $P^0$ to $U^0$ extends to a morphism $q:P\rightarrow U$.  Consider the 
diagram
\[
\xymatrix{
G\times P\ar[r]^{id\times q}\ar[d] & G\times U\ar[d]\\
P\ar[r]^q & U
}
\]
where the vertical arrows are the action maps.  Precomposing either of the two maps in the diagram 
from $G\times P$ to $U$ by the inclusion $G\times P^0\hookrightarrow G\times P$ yields the same morphism.  That is, the two 
maps from $G\times P$ to $U$ are both extensions of the same map from $G\times P^0$ to the affine scheme $U$, and hence are 
equal.  This shows that $q$ is $G$-equivariant, and therefore yields a map $g:V\rightarrow\mathcal{Y}$ such that $fg=h$.
\end{proof}
\begin{proof}[Proof of Theorem \ref{prop:cslrs}]
We begin with the following observation.  Suppose $U$ is smooth and affine over $S$ with a faithful action of a finite 
linearly reductive group scheme $G$ over $S$.  Let $y$ be a closed point of $U$ mapping to $x\in U/G$.  After making 
the \'etale base change $\Spec k(y)\rightarrow S$, we can assume $y$ is a $k$-rational point.  Let $G_y$ be 
the stabilizer subgroup scheme of $G$ fixing $y$.  Since
\[
U/G_y\longrightarrow U/G
\]
is \'etale at $y$, replacing $U/G$ by an \'etale cover, we can further assume that $G$ fixes $y$.  Then by 
Corollary \ref{cor:newcst}, we can assume $G$ has no pseudo-reflections at $y$, and hence, Theorem \ref{thm:torsor} 
shows that after shrinking $U/G$ about $x$, we can assume that the base change of $U$ to the smooth locus of 
$U/G$ is a $G$-torsor.\\
\\
We now turn to the proof.  Since $M$ has linearly reductive singularities, 
there is an \'etale cover 
$\{U_i/G_i\rightarrow M\}$, where $U_i$ is smooth and affine over $S$ and $G_i$ is a finite linearly 
reductive group scheme over $S$ which acts faithfully on $U_i$.  By the above discussion, replacing 
this \'etale cover by a finer \'etale cover if necessary, we can assume that the base change of $U_i$ to the 
smooth locus of $U_i/G_i$ is a $G_i$-torsor.  Let $M_i=U_i/G_i$ and $\mf{X}_i=[U_i/G_i]$.  We see 
that the $\mf{X}_i$ are locally the desired stacks, so we need only glue the $\mf{X}_i$.  Let 
$M_{ij}=M_i\times_M M_j$ and let $V_i\rightarrow \mf{X}_i$ be a smooth cover.  Since $M_{ij}$ is the 
coarse space of both $\mf{X}_i\times_{M_i}M_{ij}$ and $\mf{X}_j\times_{M_j}M_{ij}$, and since coarse space maps 
are dominant and codimension-preserving, Lemma \ref{l:univprop} shows that there is a unique isomorphism of 
$\mf{X}_i\times_{M_i}M_{ij}$ and $\mf{X}_j\times_{M_j}M_{ij}$.  Identifying these two stacks via this isomorphism, let 
$I_{ij}$ be the fiber product over the stack of $V_i\times_{M_i}M_{ij}$ and $V_j\times_{M_j}M_{ij}$.  We see then that we 
have a morphism $I_{ij}\rightarrow U_i\times_M U_j$.  This yields a groupoid
\[
\coprod I_{ij}\longrightarrow \coprod U_i\times_M U_j,
\]
which defines our desired glued stack $\mf{X}$.  Note that $\mf{X}$ is smooth and tame by \cite[Thm 3.2]{tame}.
\end{proof}

\end{document}